\documentclass[reqno,11pt]{amsart}
\usepackage{amsmath, latexsym, amsfonts, amssymb, amsthm, amscd}
\usepackage{epsfig}
\usepackage{graphics,epsf,psfrag}

\numberwithin{equation}{section}

\newtheorem{theorem}{Theorem}[section]
\newtheorem{lemma}[theorem]{Lemma}

\newcommand{\cA}{{\ensuremath{\mathcal A}} }

\newcommand{\cN}{{\ensuremath{\mathcal N}} }

\newcommand{\cS}{{\ensuremath{\mathcal S}} }
\newcommand{\cT}{{\ensuremath{\mathcal T}} }

\newcommand{\ga}{\alpha}

\newcommand{\gd}{\delta}

\newcommand{\gs}{\sigma}

\renewcommand{\tilde}{\widetilde}          
\DeclareMathSymbol{\leqslant}{\mathalpha}{AMSa}{"36} 
\DeclareMathSymbol{\geqslant}{\mathalpha}{AMSa}{"3E} 
\DeclareMathSymbol{\eset}{\mathalpha}{AMSb}{"3F}     

\DeclareMathOperator*{\union}{\bigcup}       
\newcommand{\suptwo}[2]{\sup_{\substack{#1 \\ #2}}} 

\newcommand{\R}{\mathbb{R}}

\newcommand{\N}{\mathbb{N}}

\def\bP{\ensuremath{\bs{\mathrm{P}}}} 
\def\bQ{\ensuremath{\bs{\mathrm{Q}}}}
\def\bE{\ensuremath{\bs{\mathrm{E}}}}

\newcommand{\ind}{\bs{1}}

\def\bs{\boldsymbol}

\newcommand{\rf}{\mathrm f}

\title{Tightness conditions for polymer measures}

\author{Francesco Caravenna}

\address{Dipartimento di Matematica Pura e Applicata,
Universit\`a degli Studi di Padova, via \mbox{Trieste} 63, 35121 Padova,
Italy } \email{francesco.caravenna\@@math.unipd.it}

\author{Giambattista Giacomin}

\address{Laboratoire de Probabilit{\'e}s et Mod\`eles Al\'eatoires  (CNRS U.M.R. 7599) and  Universit{\'e} Paris 7 -- Denis Diderot, U.F.R. Mathematiques, Case 7012, 2 place Jussieu, 75251 Paris cedex 05, France
} \email{giacomin\@@math.jussieu.fr}

\author{Lorenzo Zambotti}
\address{Laboratoire de Probabilit{\'e}s et Mod\`eles Al\'eatoires (CNRS U.M.R. 7599) and  Universit{\'e} Paris 6
-- Pierre et Marie Curie, U.F.R. Mathematiques, Case 188, 4 place Jussieu, 75252 Paris cedex 05, France
} \email{zambotti\@@ccr.jussieu.fr}

\date{\today}

\begin{document}

\begin{abstract}
We give sufficient conditions for tightness in the space $C([0,1])$ for
sequences of probability measures which enjoy a suitable decoupling between
zero level set and excursions. Applications
of our results are given in the context of (homogeneous, periodic and disordered)
random walk models for polymers and interfaces.
\\
\\
2000 \textit{Mathematics Subject Classification: 60F17, 60K35, 82B41}
\\
\\
\noindent\textit{Keywords: Tightness, Invariance Principle, Scaling Limit,
Random Walk, Polymer Model, Pinning Model, Wetting Model.}
\end{abstract}

\maketitle

\section{Introduction}

In this note we want to prove tightness under diffusive rescaling for a sequence
of processes $(\bP_N)$ with the following property: conditionally on the zero level
set, the excursions of the process between consecutive zeros are independent
and each excursion is distributed according to the same fixed law (corresponding
to the given excursion length).

\smallskip

Let us be more precise. We first need three main ingredients:
\begin{itemize}
\item the {\sl zero level set law} $p_N$ is, for each $N\in\N$,
a probability measure on the subsets of $\{1,\ldots,N\}$;
\item \rule{0pt}{1.1em}the {\sl bulk excursion law} $P_t$ is,
for each $t\in\N$, a probability measure on $\R^t$ such that
\[
P_t\left( y\in \R^t: \ y_1>0,\, \ldots, y_{t-1}>0,\, y_t=0 \right) =1\,;
\]
\item \rule{0pt}{1.1em}the {\sl final excursion law} $P_t^\rf$ is,
for each $t\in\N$, a probability measure on $\R^t$ such that
\[
P_t^\rf\left( y\in \R^t: \ y_1>0,\, \ldots,\, y_t>0 \right) =1\,.
\]
\end{itemize}
We also set for convenience
\[
\cT_N := \union_{k = 1}^N \left\{(T_1,\ldots,T_k): T_i\in\N, \
0=:T_0<T_1<\cdots<T_k\leq N \right\}\,.
\]
Then, for each $N\in\N$, we introduce the measure $\bP_N$ on $\R^N$
defined for $B_1,\ldots,B_N$ Borel sets in $\R$ by the relation:
\begin{equation}\label{eq:crucial}
\begin{split}
& \bP_N(B_1\times \ldots\times B_N) \, = \, \sum_{(T_i)\in\cT_N}
p_N(\{T_1,T_2,\ldots,T_k\}) \ \cdot
\\
& \cdot \left\{\prod_{i=1}^k
P_{T_i-T_{i-1}}\left(B_{T_{i-1}+1}\times \ldots\times
B_{T_i-T_{i-1}}\right)\right\} \cdot
P_{N-T_k}^\rf\left(B_{T_k+1}\times \ldots\times B_N\right).
\end{split}
\end{equation}
This means that under $\bP_N$:
\begin{itemize}
\item the zero level set, defined for $y\in\R^N$ by
$\cA = \cA(y):=\{\ell=1,\ldots,N: y_\ell=0\}$,
is distributed according to the law $p_N$;
\item
\rule{0pt}{1.1em}conditionally on $\cA=\{T_1,\ldots,T_k\}$, with $(T_i)\in\cT_N$, the
family of excursions $\{e_i := (y_j,j=T_{i-1}+1,\ldots,T_i)\}_{i=1,\ldots,k+1}$ is
independent;
\item
\rule{0pt}{1.1em}conditionally on $\cA=\{T_1,\ldots,T_k\}$, for all $i=1,\ldots,k$,
$e_i$ has law $P_{T_i-T_{i-1}}$; if $T_k<N$
then $e_{k+1}$ has law $P_{N-T_k}^\rf$.
\end{itemize}

We are interested in the diffusive rescaling of $\bP_N$, that we call $\bQ_N$.
More precisely, let us define the map $X^N: {\mathbb R}^N \mapsto C([0,1])$:
\[
X^N_t(y) \, = \, \frac{y_{\lfloor Nt\rfloor}}{N^{1/2}} + (Nt-\lfloor Nt\rfloor) \,
\frac {y_{\lfloor Nt\rfloor+1}-y_{\lfloor Nt\rfloor}}{N^{1/2}}, \qquad t\in[0,1],
\]
where $\lfloor r \rfloor$ denotes the integer part of $r\in\R_+$
and $y_0:=0$. Notice that $X^N_t(y)$ is nothing but the linear
interpolation of~$\{y_{\lfloor Nt \rfloor}/\sqrt{N}\}_{t \in
\frac{\N}{N} \cap [0,1]}$. Then we set
$$\bQ_N := \bP_N \circ (X^N)^{-1}\,.$$
Analogously, we denote by $Q_N$ and $Q_N^\rf$ the diffusive rescaling of
$P_N$ and $P_N^\rf$ respectively.

\smallskip

Our aim is to give sufficient conditions for tightness of $(\bQ_N)$ in $C([0,1])$:
in the next section we state and prove our main result. 
Applications to random walk models for polymers and interface are
discussed in Section~\ref{sec:polymer}.

\medskip

\section{Main result}

\smallskip
\begin{theorem}\label{main}
Suppose that:
\begin{itemize}
\item the sequences $(Q_N)$ and $(Q_N^\rf)$ are tight in $C([0,1])$;
\item the following relation holds true:
\begin{equation}\label{defK}
\lim_{a\to\infty} C(a)=0 \qquad \text{where} \
C(a) := \sup_n \ E_n \bigg( \Big( \max_{0 \le i \le n}
\frac{y_i^2}{n} \Big) \;
    \ind_{\big\{\max_{0 \le i \le n} y_i^2/n \,>\, a \big\}} \bigg)\,.
\end{equation}
\end{itemize}
Then the sequence $(\bQ_N)$ is tight in $C([0,1])$.
\end{theorem}
\smallskip

We stress that we make no hypothesis on the law $(p_N)$.
Before proving the theorem, we
introduce for $\gd > 0$ the {\sl continuity modulus} $\Gamma(\gd)$,
i.e. the real-valued functional defined for $x \in C([0,1])$ by
\begin{equation*}
    \Gamma(\gd) \;=\; \Gamma(\gd)[x] \;:=\; \suptwo{s,t \in [0,1]}{|t-s|\le \gd} |x_t - x_s| \,.
\end{equation*}
We are going to check the standard necessary and sufficient condition for tightness on $C([0,1])$
(Theorems of Prohorov and Ascoli-Arzel\`a): for every $\gamma > 0$
\begin{equation}\label{eq:tigntness0}
    \lim_{\gd \to 0} \; \sup_{N \in \N} \; \bQ_N \big( \Gamma(\gd) > \gamma \big) \;=\; 0\,.
\end{equation}
It is actually convenient to work with a modified continuity modulus $\tilde\Gamma(\gd)$:
for $x \in C([0,1])$ and $s,t \in [0,1]$ we set $s \sim_x t$ iff $x_u \ne 0$ for every $u \in (s,t)$,
i.e. iff $s$ and $t$ belong to the same excursion of $x$, and we define
\begin{equation*}
    \tilde \Gamma(\gd) \;=\; \tilde \Gamma(\gd)[x] \;:=\;
    \suptwo{s,t \in [0,1],\, s \sim_x t}{|t-s|\le \gd} |x_t - x_s| \,.
\end{equation*}
Clearly $\tilde \Gamma(\gd) \le \Gamma(\gd) \le 2 \tilde \Gamma(\gd)$, therefore it suffices to prove
that for every $\gamma > 0$
\begin{equation} \label{eq:tightness}
    \lim_{\gd \to 0} \; \sup_{N \in \N} \; \bQ_N \left( \tilde\Gamma(\gd) > \gamma \right) \;=\; 0\,.
\end{equation}

\medskip
\noindent
{\bf Proof of Theorem \ref{main}.}
The path we follow is rather general. The crucial property that we
exploit is the independence of the excursions conditionally on the
zero level set $\cA$. Setting $\cN_N :=
\union_{k=0}^\infty \{(t_1,\ldots,t_{k+1}): t_i\in\N, \
\sum_{i=1}^{k+1} t_i= N\}$, by \eqref{eq:crucial} we can write
\begin{equation}\label{eq:basic}
\begin{split}
    \bQ_N\left( \tilde \Gamma(\gd) \le \gamma \right) =
    \sum_{(t_i)\in\cN_N}
    \Bigg\{ & \prod_{\ell = 1}^k Q_{t_\ell} \Big( \Gamma\big({\textstyle \frac{N}{t_\ell}} \gd\big)
    \le \gamma {\textstyle \sqrt{ \frac{N}{t_\ell}} } \Big) \Bigg\}
    \cdot Q_{t_{k+1}}^\rf \Big(\, \Gamma\big({\textstyle \frac{N}{t_{k+1}}} \gd\big)
    \le \gamma {\textstyle \sqrt{ \frac{N}{t_{k+1}}} } \,\Big)
    \\
    & \cdot \, p_N \big( \{t_1,t_1+t_2,\ldots,t_1+\cdots+t_{k+1}\}\big)\,.
\end{split}
\end{equation}
Next we perform a very drastic bound: we set
\begin{equation}\label{eq:def_f}
    f_\gamma(\gd) \;:=\; \inf_{N\in\N,\, (t_i)\in\cN_N}
    \     \prod_{\ell = 1}^k Q_{t_\ell} \Big(\, \Gamma\big({\textstyle \frac{N}{t_\ell}} \gd\big)
    \le \gamma {\textstyle \sqrt{ \frac{N}{t_\ell}} } \,\Big)\,,
\end{equation}
\[
g_\gamma(\delta) \; := \; \inf_{N\in\N, \, 1\leq t\leq N} \
Q_{t}^\rf \left(\, \Gamma\big({\textstyle \frac{N}{t}} \gd\big) \le
\gamma {\textstyle \sqrt{ \frac{N}{t}} }\right).
\]
By \eqref{eq:basic} we have $\bQ_N( \tilde \Gamma(\gd) \le \gamma )
\, \ge \, f_\gamma(\gd)\cdot g_\gamma(\gd)$. Therefore if we show
that $f_\gamma(\gd)g_\gamma(\gd) \to 1$ as $\gd \to 0$, for any
fixed $\gamma > 0$, equation \eqref{eq:tightness} follows and the
proof is completed.

\medskip
We start by proving that $\liminf_{\delta\searrow 0}
f_\gamma(\delta)=1$ for all $\gamma>0$. We introduce an auxiliary (small) parameter
$\eta
> 0$ and we define the set
\begin{equation*}
    \cS^\eta \;=\; \cS^\eta (N,(t_i)) \;:=\;
    \big\{ \ell \in \{1,\ldots,k\}:\; t_\ell > \eta N \big\}\,,
    \qquad N\in\N, \ (t_i)\in\cN_N.
\end{equation*}
Notice that $\Gamma(\cdot)$ is non-decreasing and that
we have trivially $\Gamma(\gd')[x] \le 2 \max_{t \in [0,1]} |x_t|$.
Splitting the product in the r.h.s. of \eqref{eq:def_f} and using these observations, we obtain
\begin{equation}\label{eq:split}
\begin{split}
    \prod_{\ell = 1}^{k} Q_{t_\ell} \Big(\, \Gamma\big({\textstyle \frac{N}{t_\ell}} \gd\big)
    \le \gamma {\textstyle \sqrt{ \frac{N}{t_\ell}} } \,\Big) \;\ge\;
    \Bigg\{ & \prod_{\ell \in \cS^\eta} Q_{t_\ell} \Big(\, \Gamma\big({\textstyle \frac{\gd}{\eta}} \big)
    \le \gamma \,\Big) \Bigg\} \\
    & \cdot \Bigg\{ \prod_{\ell \in \{1, \ldots, k\} \setminus \cS^\eta}
    P_{t_\ell} \Big(\, \max_{0 \le i \le t_\ell} |y_i|
    \le {\textstyle \frac\gamma 2} \sqrt{N} \,\Big) \Bigg\} \,.
\end{split}
\end{equation}

Suppose now that we can prove the following:
\begin{equation}\label{eq:claimm}
\forall \ \gamma>0: \qquad \lim_{\eta\searrow 0} \
\inf_{N,(t_i)\in\cN_N} \ \prod_{\ell \in \{1, \ldots, k\} \setminus
\cS^\eta}
    P_{t_\ell} \Big(\, \max_{0 \le i \le t_\ell} |y_i|
    \le {\textstyle \frac\gamma 2} \sqrt{N} \,\Big) = 1.
\end{equation}
In other words, for any fixed $\gamma> 0$, the parameter $\eta$ can
be chosen in order to make the second term in the r.h.s. of \eqref{eq:split}
as close to $1$ as we wish, {\sl uniformly in $N$ and $(t_i)$}.
If \eqref{eq:claimm} is proven, then we can fix $\eta>0$ and it is
easy to see that one can choose $\gd$ in order to make also the
first term in the r.h.s. of \eqref{eq:split} as close to $1$ as we wish, {\sl uniformly
in $N$ and $(t_i)$}: this is just because the number of factors in
the product (i.e. the cardinality of the set $\cS^\eta$) is bounded
by construction by $1/\eta < \infty$ and because by hypothesis the
sequences $(Q_N)$ and $(Q_N^\rf)$ are tight in $C([0,1])$ (we recall
that $t_\ell \ge \eta N$ for $\ell \in \cS^\eta$). This shows that
indeed $f_\gamma(\gd) \to 1$ as $\gd \to 0$, for any fixed $\gamma >
0$, completing the proof.

Therefore we are left with proving \eqref{eq:claimm}: we have to
show that for any fixed $\gamma$ and $\ga > 0$ we can choose the
parameter $\eta$ such that
\begin{equation} \label{eq:claim0}
    \inf_{N\in\N,\, (t_i)\in\cN_N}
    \Bigg( \prod_{\ell \in \{1, \ldots, k\} \setminus \cS^\eta}
    P_{ t_\ell} \Big(\, \max_{0 \le i \le t_\ell} |y_i|
    \le {\textstyle \frac\gamma 2} \sqrt{N} \,\Big) \Bigg) \;\ge\; 1 - \ga \,.
\end{equation}
By (a somewhat enhanced) Chebychev inequality we have
\begin{equation*}
    P_{n} \bigg(\, \max_{0 \le i \le n} \frac{|y_i|}{\sqrt n} > a \,\bigg) \;\le\;
    \frac{c_n(a)}{a^2}\,, \qquad c_n(a) \;:=\;
    \bE_{n}\bigg( \Big( \max_{0 \le i \le n} \frac{y_i^2}{n} \Big) \;
    \ind_{\big\{\max_{0 \le i \le n} y_i^2/n \,>\, a \big\}} \bigg)\,.
\end{equation*}
Since $c_n(\cdot)$ is non-increasing and $t_\ell \le \eta N$ for $\ell \not \in \cS^\eta$,
we obtain
\begin{equation} \label{eq:step1}
    \prod_{\ell \in \{1, \ldots, k\} \setminus \cS^\eta}
    P_{ t_\ell} \Big(\, \max_{0 \le i \le t_\ell} |y_i|
    \le {\textstyle \frac\gamma 2} \sqrt{N} \,\Big) \;\ge\;
    \prod_{\ell \in \{1, \ldots, k\} \setminus \cS^\eta} \bigg( 1 -
    c_{t_\ell} \Big( {\textstyle \frac{\gamma}{2\sqrt\eta} } \Big)
    \frac{4}{\gamma^2} \, \frac{t_\ell}{N} \bigg)\,.
\end{equation}
Notice now that $C(a)=\sup_{n\in\N} c_n(a)$ by the definition
\eqref{defK} of $C$. Choosing $\eta$ sufficiently small, so that
$c_{t_\ell} \big( \frac{\gamma}{2\sqrt\eta} \big) \frac{4}{\gamma^2}
\frac{t_\ell}{N} \,\le\, C \big( \frac{\gamma}{2\sqrt\eta} \big)
\frac{4}{\gamma^2} \eta \,\le\, \frac 12 \,$, and observing that
$1-x \ge \exp(-2x)$ for $0 \le x \le \frac 12$, we can finally bound
\eqref{eq:step1} by
\begin{equation} \label{eq:step2}
\begin{split}
    \prod_{\ell \in \{1, \ldots, k\} \setminus \cS^\eta} &
    P_{ t_\ell} \Big(\, \max_{0 \le i \le t_\ell} |y_i|
    \le {\textstyle \frac\gamma 2} \sqrt{N} \,\Big) \;\ge\;
    \prod_{\ell \in \{1, \ldots, k\} \setminus \cS^\eta} \exp \bigg( - 2
    C \Big( {\textstyle \frac{\gamma}{2\sqrt\eta} } \Big)
    \frac{4}{\gamma^2} \, \frac{t_\ell}{N} \bigg)\\
    & \;=\; \exp \bigg( - 2 C \Big( {\textstyle \frac{\gamma}{2\sqrt\eta} } \Big)
    \frac{4}{\gamma^2} \sum_{\ell \in \{1, \ldots, k\} \setminus \cS^\eta}
    \frac{t_\ell}{N} \bigg) \;\ge\;
    \exp \bigg( - \frac{8}{\gamma^2} \, C \Big( {\textstyle \frac{\gamma}{2\sqrt\eta} } \Big) \bigg) \,,
\end{split}
\end{equation}
and equation \eqref{eq:claim0} follows from the hypothesis
$\lim_{a\to\infty} C(a)=0$.

\medskip
It remains to prove that
$\liminf_{\delta\searrow 0} g_\gamma(\delta)=1$, for all $\gamma>0$. Let
$\theta\in(0,1)$. Then
\[
g^{(1)}_\gamma(\delta) := \inf_{N\in\N, \, \theta N\leq t\leq N} \
Q_{t}^\rf \left(\, \Gamma\big({\textstyle \frac{N}{t}} \gd\big) \le
\gamma {\textstyle \sqrt{ \frac{N}{t}} }\right) \geq \inf_t
Q_{t}^\rf \left(\, \Gamma\big({\textstyle \frac{\delta}{\theta}}
\big) \le \gamma \right)
\]
so that, by tightness of $(Q_N^\rf)$ in $C([0,1])$, we have
$\liminf_{\delta\searrow 0} g_\gamma^{(1)}(\delta)=1$ for all
$\gamma>0$ and $\theta\in(0,1)$. Now:
\[
g^{(2)}_\gamma(\delta) := \inf_{N\in\N, \, 1\leq t< \theta N} \
Q_{t}^\rf \left(\, \Gamma\big({\textstyle \frac{N}{t}} \gd\big) \le
\gamma {\textstyle \sqrt{ \frac{N}{t}} }\right) \geq \inf_t
Q_{t}^\rf \left(\, \sup_{s\in[0,1]} \, |x_s|  \le
\frac\gamma{2\sqrt\theta} \right),
\]
and again by tightness of $(Q_N^\rf)$ in $C([0,1])$, we have that
for any $\alpha\in(0,1)$ we can find $\theta\in(0,1)$ such that
$\liminf_{\delta\searrow 0} g_\gamma^{(2)}(\delta)\geq 1-\alpha$. It
follows that $\liminf_{\delta\searrow 0} g_\gamma(\delta)\geq
1-\alpha$ for all $\alpha\in(0,1)$, and the proof is completed.\qed

\medskip

\section{Application to polymer measures}
\label{sec:polymer}

One direct application of Theorem \ref{main}, and the main motivation
of this note, is in the context of $(1+1)$--dimensional
random walk models for polymer chains
and interfaces. We look in particular at the {\sl copolymer near a selective
interface model}, both in the disordered \cite{cf:BdH,cf:BG2} and in the
periodic \cite{cf:BG,cf:CGZ1} setting, but also at the {\sl interface wetting models}
considered in \cite{cf:DGZ,cf:CGZ2} and at {\sl pinning models} based on random
walks, described e.g. in \cite{cf:Gia} (to which we refer for a detailed overview
on all these models).

Notice that, for the purpose of proving tightness in $C[0,1]$, one can safely focus
on the absolute value of the process. With this observation in mind, we have
the basic fact
that all the above mentioned models satisfy equation \eqref{eq:crucial},
for suitable choices of the laws $p_N$, $P_N$ and $P_N^\rf$. More precisely,
in all these cases we have that for every Borel set $B \in \R^t$
\begin{equation} \label{eq:twolaws}
\begin{split}
	P_t (B) &= \bP \big( (S_1, \ldots, S_t) \in B \,\big|\, S_1>0,\ldots,S_{t-1}>0, S_t=0 \big)\\
	P_t^\rf(B) &= \bP \big( (S_1, \ldots, S_t) \in B \,\big|\, S_1>0, \ldots, S_t>0 \big)\,,
\end{split}
\end{equation}
where $(\{S_n\}_{n\ge 0}, \bP)$ is a real random walk.
(In particular, in the non-homogeneous cases,
all the dependence on the environment is contained in the law $p_N$.)
Then by Theorem~\ref{main}
the tightness under diffusive rescaling for all the above models
is reduced to showing that the
sequences $(P_N)$ and $(P_N^\rf)$ are tight and that equation \eqref{defK} holds.

We are going to check these conditions in the special instance when $(\{S_n\}_{n\ge 0}, \bP)$
is a non-trivial symmetric random walk with $S_1 \in \{-1,0,+1\}$, thus
proving tightness for the (disordered and periodic) copolymer near a selective interface model.
Notice that the law of the walk is identified by
$p:=\bP \left( S_1 = +1\right) \in (0, 1/2]$.
To lighten notations, in what follows we actually assume that $p \in (0,1/2)$.
The tightness for the sequences $(Q_N)$ and $(Q_N^\rf)$ in this case is a classical
result, cf. \cite{cf:K} and \cite{cf:B}. Therefore it remains to prove that equation \eqref{defK}
holds true, which follows as a simple consequence of the following lemma.

\smallskip
\begin{lemma}\label{le3.1} There exists a constant $C>0$ such that for all $n\in\N$ and
$a> 0$ we have:
\[
f_n(a):=\bP\left( \max_{1\leq i\leq n} \frac{S^2_i}n \geq  a \ \Big| \ S_1 > 0, \ldots, S_{n-1} > 0, S_n=0
\right) \leq C \, \frac 1{1+ a^2}.
\]
\end{lemma}
\smallskip
\noindent {\bf Proof}.
Set $\tau:=\inf\{n>0: S_n^2\geq na\}$ and
$T:=\inf\{n>0: S_n\leq 0\}$. Then for $n\geq 1$:
\[
f_n(a) = \bP\left( \tau\leq n \ \Big| \ T= n
\right) = \frac{\bP\left( \tau\leq n, \ T= n\right)}{\bP\left( T=n \right)}
\]
By the reflection principle, the denominator is equal to:
\[
\bP\left( T=n \right) = p^2 \left[ \bP(S_{n-2}=0)-\bP(S_{n-2}=2) \right].
\]
By symmetry and by the strong Markov property, we can estimate the numerator:
\[
\begin{split}
\bP\left( \tau\leq n, \ T= n\right) \ \leq & \ 2 \,
\bP\left( \tau\leq n/2, \ T= n \right)
\\ = & \ 2 \sum_{j=1}^{n/2} \bP\left(\tau=j\leq n/2< T\right)  \
\bP( S_{n-j}= a_n, \ T>n-j)
\end{split}
\]
where $a_n:=\lfloor \sqrt{an}\rfloor$ is the integer part of $\sqrt{an}$.
Again by the reflection principle:
\[
\bP( S_{n-j}= a_n, \ T>n-j) = p
\left[
\bP \left( S_{n-j-1}=a_n-1 \right) -
\bP  \left( S_{n-j-1}=a_n+1 \right)
\right].
\]
Th. 16 in Ch. VII of \cite{petrov}
says that when $k\to\infty$, uniformly in $b\in \N/\sqrt{k}$:
\[
\left( 1+ \left\vert \frac b\gs \right\vert^4\right)
\left( \gs \sqrt{k} \bP \left(S_k= b \sqrt{k}\right)
- \frac 1{\sqrt{2\pi}} \exp(-b^2/(2\gs))-
\sum_{\nu=1}^2 \frac{q_\nu (b/\gs)}{k^{\nu/2}}
\right) = o\left( k^{-1}\right),
\]
where $\gs =\sqrt{2p}$.
By using the last two formulas we obtain that
there exist positive constants $c_1$ and  $c_2$ such that for all $j\leq n/2$
\begin{equation}
\bP \left( T>n-j, S_{n-j}= a_n\right) \, \le \,
\frac{c_1}n \left( \exp\left(-c_2 a\right)+  \frac 1{1+ a^2}\right), \qquad
\forall \ a>0.
\end{equation}
Now:
\[
\sum_{j=1}^{n/2} \bP\left(\tau=j\leq n/2< T\right) \leq
\bP\left(\tau< T\right) = \frac 1{1+a_n} \leq \frac 1{\sqrt{an}},
\]
where the equality can be proved with a martingale argument.
Similarly, we have that
\begin{equation} \label{eq:as_K}
    \exists \lim_{n\to\infty} n^{3/2} \, \bP(T=n) \;=:\; c_K \in (0,\infty)\,.
\end{equation}
(see also \cite[Ch. XII.7]{cf:Feller2}).
It is now easy to conclude. \qed

\end{document}